\documentclass{article}
\pdfoutput=1
\usepackage{graphicx}
\usepackage{subcaption}
\usepackage[includefoot,margin=1in]{geometry} 

\usepackage{lscape}
\usepackage[english]{babel} 

\usepackage{tikz}
\usetikzlibrary{arrows.meta}
\usetikzlibrary{plotmarks}
\usetikzlibrary{shapes}
\usetikzlibrary{arrows}
\usepackage{pgfplots}
\pgfplotsset{compat=newest}
 
\usepackage{array} 
\usepackage{footnote}
 
\usepackage{amssymb}
\usepackage{amsmath}
\usepackage{amsfonts}
\usepackage{bm}

\usepackage{url}
\usepackage{hyperref}
\hypersetup{pdfstartview={XYZ null null 1.00}}

\newcommand{\beq}{\begin{equation}}
\newcommand{\eeq}{\end{equation}}

\newcommand{\beqr}{\begin{eqnarray}}
\newcommand{\eeqr}{\end{eqnarray}}
\def\bal#1\eal{\begin{align}#1\end{align}}
\def\bat#1#2\eat{\begin{alignat}{#1}#2\end{alignat}}

\begin{document}

\title{Folding $\pi$}
\author{Michael Assis}
\date{Department of Obstetrics, Gynaecology and Newborn Health, The University of Melbourne\\[2ex]\today}
\maketitle
\begin{abstract}
It is well known that the set of origami constructible numbers is larger than the classical straight-edge and compass constructible numbers. However, the Huzita-Justin-Hatori origami constructible numbers remain algebraic so that the transcendental number $\pi$ can only be approximated using a finite number of straight line folds. Using these methods we give a convergent sequence for folding $\pi$ as well as other methods to approximate $\pi$. Folding along curved creases, however, allows for the construction of transcendental numbers. We here give a method to construct $\pi$ exactly by folding along a parabola, and we discuss generalizations for folding other transcendental numbers such as $\Gamma(1/4)$. 
\end{abstract}

\section{Introduction}
The number $\pi$ is a transcendental number, which means that it isn't the solution of any algebraic equation~\cite{klein1897}. From the perspective of classical constructibility, this means that $\pi$ cannot be constructed, since straight-edge and compass can only construct numbers which are solutions of quadratic equations or nested square roots of such numbers~\cite{klein1897}. Origami can construct a larger class of numbers using the Huzita-Justin-Hatori origami operations~\cite{alperin2009}, which allows the construction of numbers derived from cubic equations. However, just like with straight-edge and compass, straight-line creases cannot construct $\pi$ exactly. To go further, curved creases are necessary. 

We begin by first considering methods to construct $\pi$ approximately using straight-line creases. Just as with the numerical algorithms that can keep running to give ever more precise values, we begin by giving a convergent sequence whose approximation is limited only by the physical limitations of folding ever finer angles with precision. On the other hand, we can also consider particular cut-off points of precision to $\pi$ and attempt to fold these approximations. We give three such methods, first folding $\pi$'s continued fraction convergents using Fujimoto's approximation method, then folding trigonometric approximants using an angular version of Fujimoto's approximation method, and finally folding cubic approximants to $\pi$. 
 
All of the straight-line crease methods only approximate $\pi$. We therefore expand into allowing the folding of curved creases. Of course,  folding along a curved crease no longer follows the Huzita-Justin axioms of straight-line crease origami operations, since they cannot construct curved creases exactly. However, by folding along a parabola, we show that we can construct $\pi$ exactly. Finally, we end by discussing constructions for other transcendental numbers and how to view the curved crease construction as a  convergent sequence of straight-line origami operations.

\section{Convergent folding sequence}
In 2007, Osler gave a geometric translation of Vieta's product formula for $\pi$~\cite{osler2007}, 
\beq
\frac{2}{\pi} = \sqrt{\frac{1}{2}}\sqrt{\frac{1}{2}+\frac{1}{2}\sqrt{\frac{1}{2}}}\sqrt{\frac{1}{2}+\frac{1}{2}\sqrt{\frac{1}{2}+\frac{1}{2}\sqrt{\frac{1}{2}}}}\cdots
\eeq
which can be rewritten as 
\beq
\frac{\pi}{2} = \prod_{n=1}^{\infty} \sec\left(\frac{\pi}{2^{n+1}} \right)
\eeq
Osler's geometric interpretation of this equation can be readily adapted to an origami folding sequence. We show such an adaptation in Figure~\ref{figure1} for a square sheet of paper, such as typical origami paper, although using rectangles such as Letter of A4 paper can also be used.

In practice, only 6 to 8 iterations can be made due to the difficulty of folding small angles accurately. With 6 iterations, the the procedure will give 3 accurate digits, $3.14127\ldots$, while with 8 iterations, the procedure will give 4 accurate digits, $3.141572\ldots$. If a sheet of paper is cut to exactly $4\times4$~cm, then the $\pi$ approximant can be measured directly in cm from the last iteration, although in practice this size paper is too small. In general, the approximant to $\pi$ is found by measuring the last iteration length and dividing by the length of the 1/4 horizontal mark, shown in the last step.  

Due to folding errors, crease size, and measuring limitations, perhaps only 2 decimals can be found this way in practice, although with a very large sheet 3 or 4 digits might be found. Photographic backdrop paper is inexpensive, colorful, and typically comes in rolls of width around 9 ft (3 m), which would be suitable for a large demonstration.

\begin{figure*}[htpb]
\begin{center}
\begin{subfigure}{0.45\textwidth}
\centering
\begin{tikzpicture}[scale=1.4]
\draw[line width = 1pt] (0,0) rectangle (4,4);
\draw[loosely dashed,line width = 1pt] (0,0) -- (4,4);
\node (a) at (3.9,0) {};
\node (b) at (0,3.9) {};
\draw[<->] (a)  to [out=90,in=0, looseness=1.2] (b);
\end{tikzpicture}
\caption{Fold from corner to corner\label{fig1:step1}}
\end{subfigure}
\hspace{0.08\textwidth}
\begin{subfigure}{0.45\textwidth}
\centering
\begin{tikzpicture}[scale=1.4]
\draw[line width = 1pt] (0,0) rectangle (4,4);
\draw[line width = 0.5pt] (0,0) -- (4,4);
\draw[loosely dashed,line width = 1pt] (0,4) -- (4,0);
\node (a) at (.1,.1) {};
\node (b) at (3.9,3.9) {};
\draw[<->] (a)  to [out=90,in=180, looseness=1.2] (b);
\end{tikzpicture}
\caption{Fold from corner to corner \label{fig1:step2}}
\end{subfigure}
\\[0.2in]
\begin{subfigure}{0.45\textwidth}
\centering
\begin{tikzpicture}[scale=1.4]
\draw[line width = 1pt] (0,0) rectangle (4,4);
\draw[line width = 0.5pt] (0,0) -- (4,4);
\draw[line width = 0.5pt] (0,4) -- (4,0);
\draw[loosely dashed,line width = 1pt] (0,0) -- (4,1.65685);
\node (a) at (1.1,1.1) {};
\node (b) at (1.5,0.0) {};
\draw[<->] (a)  to [out=315,in=90, looseness=1.2] (b);
\end{tikzpicture}
\caption{Bisect the angle \label{fig1:step3}}
\end{subfigure}
\hspace{0.08\textwidth}
\begin{subfigure}{0.45\textwidth}
\centering
\begin{tikzpicture}[scale=1.4]
\draw[line width = 1pt] (0,0) rectangle (4,4);
\draw[line width = 0.5pt] (0,0) -- (4,4);
\draw[line width = 0.5pt] (0,4) -- (4,0);
\draw[line width = 0.5pt] (0,0) -- (4,1.65685);
\draw[dashed,line width = 1pt] (2.8284,1.17157) -- (3.313708,0);
\draw[line width = 0.5pt,rotate around={22.5:(2.8284,1.17157)}] (2.8284,1.17157) rectangle (2.8284-0.2,1.17157-0.2);
\node (a) at (4,1.65685) {};
\node (b) at (1.6,0.7) {};
\draw[<->] (a)  to [out=135,in=90, looseness=1.1] (b);
\end{tikzpicture}
\caption{Fold the previous crease back on itself through the intersection \label{fig1:step4}}
\end{subfigure}
\\[0.2in]
\begin{subfigure}{0.45\textwidth}
\centering
\begin{tikzpicture}[scale=1.4]
\draw[line width = 1pt] (0,0) rectangle (4,4);
\draw[line width = 0.5pt] (0,0) -- (4,4);
\draw[line width = 0.5pt] (0,4) -- (4,0);
\draw[line width = 0.5pt] (0,0) -- (4,1.65685);
\draw[line width = 0.5pt] (2.8284,1.17157) -- (3.313708,0);
\draw[dashed,line width = 1pt] (0,0) -- (4,0.79565);
\node (a) at (1.7,0.8) {};
\node (b) at (2,0.0) {};
\draw[<->] (a)  to [out=315,in=90, looseness=1.2] (b);
\end{tikzpicture}
\caption{Bisect the angle \label{fig1:step5}}
\end{subfigure}
\hspace{0.08\textwidth}
\begin{subfigure}{0.45\textwidth}
\centering
\begin{tikzpicture}[scale=1.4]
\draw[line width = 1pt] (0,0) rectangle (4,4);
\draw[line width = 0.5pt] (0,0) -- (4,4);
\draw[line width = 0.5pt] (0,4) -- (4,0);
\draw[line width = 0.5pt] (0,0) -- (4,1.65685);
\draw[line width = 0.5pt] (2.8284,1.17157) -- (3.313708,0);
\draw[line width = 0.5pt] (0,0) -- (4,0.79565);
\draw[dashed,line width = 1pt] (3.06147,0.60896) -- (3.1826,0);
\draw[line width = 0.5pt,rotate around={11.25:(3.06147,0.60896)}] (3.06147,0.60896) rectangle (3.06147-0.2,0.60896-0.2);
\node (a) at (4,0.79565) {};
\node (b) at (2,0.4) {};
\draw[<->] (a)  to [out=135,in=90, looseness=0.6] (b);
\end{tikzpicture}
\caption{Fold the previous crease back on itself through the nearest intersection \label{fig1:step6}}
\end{subfigure}
\end{center}
\end{figure*}

\begin{figure}[htpb]
\addtocounter{figure}{-1}
\begin{center}
\begin{subfigure}{0.45\textwidth}
\addtocounter{subfigure}{6}
\centering
\begin{tikzpicture}[scale=1.4]
\draw[line width = 1pt] (0,0) rectangle (4,4);
\draw[line width = 0.5pt] (0,0) -- (4,4);
\draw[line width = 0.5pt] (0,4) -- (4,0);
\draw[line width = 0.5pt] (0,0) -- (4,1.65685);
\draw[line width = 0.5pt] (2.8284,1.17157) -- (3.313708,0);
\draw[line width = 0.5pt] (0,0) -- (4,0.79565);
\draw[line width = 0.5pt] (3.06147,0.60896) -- (3.1826,0);
\node (a) at (2,0.13) {$\bm{\ldots}$};
\end{tikzpicture}
\caption{Repeat previous two steps as much as practicable\label{fig1:step7}}
\end{subfigure}
\hspace{0.08\textwidth}
\begin{subfigure}{0.45\textwidth}
\centering
\begin{tikzpicture}[scale=1.4]
\draw[line width = 1pt] (0,0) rectangle (4,4);
\draw[line width = 0.5pt] (0,0) -- (4,4);
\draw[line width = 0.5pt] (0,4) -- (4,0);
\draw[line width = 0.5pt] (0,0) -- (4,1.65685);
\draw[line width = 0.5pt] (2.8284,1.17157) -- (3.313708,0);
\draw[line width = 0.5pt] (0,0) -- (4,0.79565);
\draw[line width = 0.5pt] (3.06147,0.60896) -- (3.1826,0);
\draw[line width = 0.5pt] (0,0) -- (4,0.393965);
\draw[line width = 0.5pt] (3.121445,0.307436) -- (3.151725,0);
\draw[line width = 0.5pt] (0,0) -- (4,0.196507);
\draw[line width = 0.5pt] (3.136548,0.154089) -- (3.144118,0);
\draw[line width = 0.5pt] (0,0) -- (4,0.0981945);
\draw[line width = 0.5pt] (3.140331,0.077091) -- (3.142224,0);
\draw[line width = 0.5pt] (0,0) -- (4,0.04909);
\draw[line width = 0.5pt] (3.14128,0.03855) -- (3.14175,0);
\draw[dashed,line width = 1pt] (1,0) -- (1,4);
\node (a) at (0,1.9) {};
\node (b) at (2,1.9) {};
\draw[<->] (a)  to [out=45,in=135, looseness=0.6] (b);
\end{tikzpicture}
\caption{Fold the left edge to the center \label{fig1:step8}}
\end{subfigure}
\\[0.2in]
\begin{subfigure}{0.9\textwidth}
\centering
\begin{tikzpicture}[scale=2]
\draw[line width = 1pt] (0,0) rectangle (4,4);
\draw[line width = 0.5pt] (0,0) -- (4,4);
\draw[line width = 0.5pt] (0,4) -- (4,0);
\draw[line width = 0.5pt] (0,0) -- (4,1.65685);
\draw[line width = 0.5pt] (2.8284,1.17157) -- (3.313708,0);
\draw[line width = 0.5pt] (0,0) -- (4,0.79565);
\draw[line width = 0.5pt] (3.06147,0.60896) -- (3.1826,0);
\draw[line width = 0.5pt] (0,0) -- (4,0.393965);
\draw[line width = 0.5pt] (3.121445,0.307436) -- (3.151725,0);
\draw[line width = 0.5pt] (0,0) -- (4,0.196507);
\draw[line width = 0.5pt] (3.136548,0.154089) -- (3.144118,0);
\draw[line width = 0.5pt] (0,0) -- (4,0.0981945);
\draw[line width = 0.5pt] (3.140331,0.077091) -- (3.142224,0);
\draw[line width = 0.5pt] (0,0) -- (4,0.04909);
\draw[line width = 0.5pt] (3.14128,0.03855) -- (3.14175,0);
\draw[line width = 0.5pt] (1,0) -- (1,4);
\draw[<->,thin] (0.05,-.3) -- (0.95,-0.3); 
\node at (0.5,-0.23) {$a$};
\draw[<->,thin] (0.05,-.5) -- (3.0915,-0.5); 
\node at (1.57,-0.4) {$b$};
\draw[thin] (0,-0.1) -- (0,-.6);
\draw[thin] (1,-.1) -- (1,-.4);
\draw[thin] (3.1415,-0.1) -- (3.1415,-.6);
\end{tikzpicture}
\caption{Measure lengths $a$ and $b$. In the limit, $b/a\to \pi$.\label{fig1:step9}}
\end{subfigure}
\end{center}
\caption{Convergent folding sequence for $\pi$.\label{figure1}}
\end{figure}
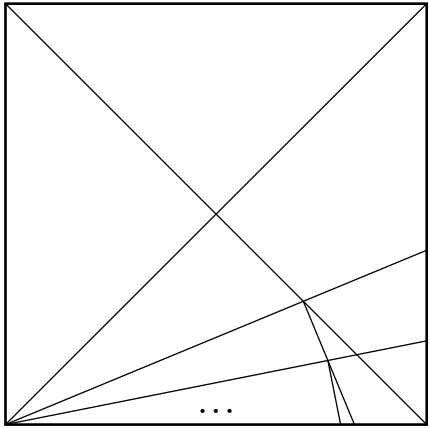
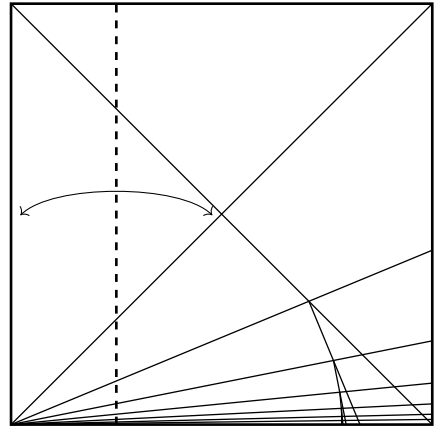
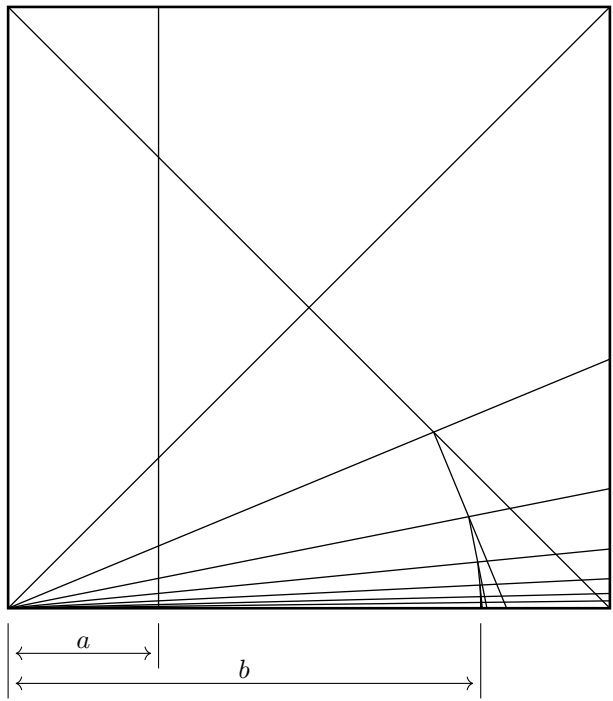

\newpage
\section{Folding \texorpdfstring{$\pi$}{Pi}'s continued fraction convergents}
\begin{figure}[htpb]
\begin{center}
\begin{subfigure}{0.3\textwidth}
\centering
\begin{tikzpicture}[scale=1]
\draw[line width = 1pt] (0,0) rectangle (5,1);
\draw[dashed, line width = 1pt] (1.2,0) -- (1.2,1);
\node at (1.2,1.3) {$\frac{1}{5}+\epsilon$};
\end{tikzpicture}
\caption{Initial guess for $\frac{1}{5}$ with error $\epsilon$. \label{fig2:step1}}
\end{subfigure}
\hspace{0.03\textwidth}
\begin{subfigure}{0.3\textwidth}
\centering
\begin{tikzpicture}[scale=1]
\draw[line width = 1pt] (0,0) rectangle (5,1);
\draw[line width = 0.5pt] (1.2,0) -- (1.2,1);
\draw[dashed, line width = 1pt] (3.1,0) -- (3.1,1);
\node at (3.1,1.3) {$\frac{3}{5}+\frac{\epsilon}{2}$};
\node (a) at (5,0.4) {};
\node (b) at (1.2,0.4) {};
\draw[<->] (a)  to [out=135,in=45, looseness=0.4] (b);
\end{tikzpicture}
\caption{Fold the right edge to the previous crease (R) \label{fig2:step2}}
\end{subfigure}
\hspace{0.03\textwidth}
\begin{subfigure}{0.3\textwidth}
\centering
\begin{tikzpicture}[scale=1]
\draw[line width = 1pt] (0,0) rectangle (5,1);
\draw[line width = 0.5pt] (1.2,0) -- (1.2,1);
\draw[line width = 0.5pt] (3.1,0) -- (3.1,1);
\draw[dashed, line width = 1pt] (4.05,0) -- (4.05,1);
\node at (4.05,1.3) {$\frac{4}{5}+\frac{\epsilon}{4}$};
\node (a) at (5,0.4) {};
\node (b) at (3.1,0.4) {};
\draw[<->] (a)  to [out=135,in=45, looseness=0.4] (b);
\end{tikzpicture}
\caption{Fold the right edge to the previous crease (R)\label{fig2:step3}}
\end{subfigure}
\\[0.2in]
\begin{subfigure}{0.3\textwidth}
\centering
\begin{tikzpicture}[scale=1]
\draw[line width = 1pt] (0,0) rectangle (5,1);
\draw[line width = 0.5pt] (1.2,0) -- (1.2,1);
\draw[line width = 0.5pt] (3.1,0) -- (3.1,1);
\draw[line width = 0.5pt] (4.05,0) -- (4.05,1);
\draw[dashed, line width = 1pt] (2.025,0) -- (2.025,1);
\node at (2.025,1.3) {$\frac{2}{5}+\frac{\epsilon}{8}$};
\node (a) at (4.05,0.4) {};
\node (b) at (0,0.4) {};
\draw[<->] (a)  to [out=135,in=45, looseness=0.4] (b);
\end{tikzpicture}
\caption{Fold the left edge to the previous crease (L) \label{fig2:step4}}
\end{subfigure}
\hspace{0.03\textwidth}
\begin{subfigure}{0.3\textwidth}
\centering
\begin{tikzpicture}[scale=1]
\draw[line width = 1pt] (0,0) rectangle (5,1);
\draw[line width = 0.5pt] (1.2,0) -- (1.2,1);
\draw[line width = 0.5pt] (3.1,0) -- (3.1,1);
\draw[line width = 0.5pt] (4.05,0) -- (4.05,1);
\draw[line width = 0.5pt] (2.025,0) -- (2.025,1);
\draw[dashed, line width = 1pt] (1.0125,0) -- (1.0125,1);
\node at (1.0125,1.3) {$\frac{1}{5}+\frac{\epsilon}{16}$};
\node (a) at (2.025,0.4) {};
\node (b) at (0,0.4) {};
\draw[<->] (a)  to [out=135,in=45, looseness=0.4] (b);
\end{tikzpicture}
\caption{Fold the left edge to the previous crease (L) \label{fig2:step5}}
\end{subfigure}
\hspace{0.03\textwidth}
\begin{subfigure}{0.3\textwidth}
\centering
\begin{tikzpicture}[scale=1]
\end{tikzpicture}
\end{subfigure}
\end{center}
\caption{Folding Fujimoto's exponentially convergent approximation for $\frac{1}{5}$ from an initial guess of $\frac{1}{5}$ with error $\epsilon$. The steps taken can be represented as RRLL. \label{figure2}}
\end{figure}
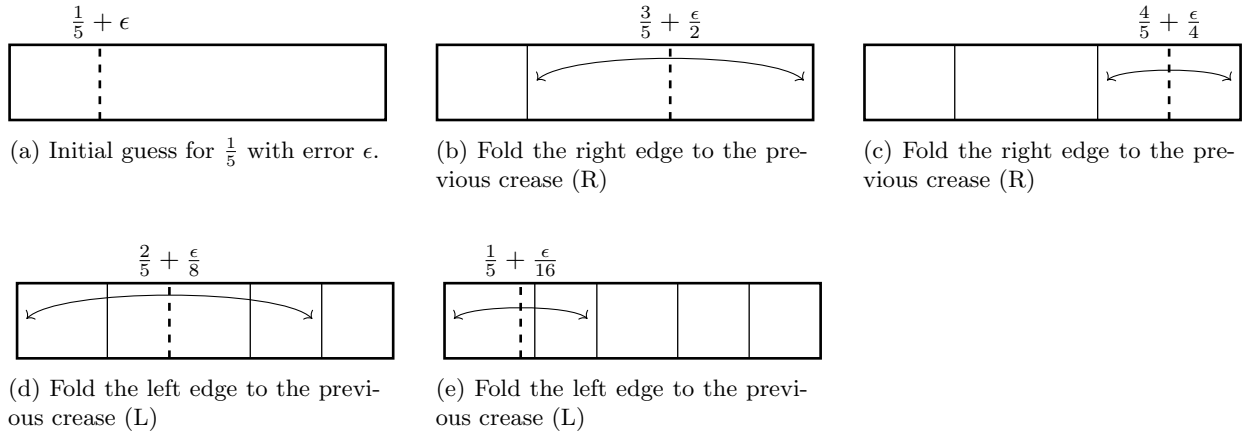

The Fujimoto approximation~\cite{fujimoto1978,fujimoto1982} is an exponentially converging algorithm to divide a strip of paper into a rational fraction with odd denominator of a length between 0 and 1. Beginning with an initial guess, a series of folds is made, such that either the left or the right edge of the paper is taken to the previously made fold; see Figure~\ref{figure2} for a demonstration involving approximating $\frac{1}{5}$. After the procedure returns to near the original guess, a new more accurate crease will have been created much closer to the desired rational fraction. The process can be further repeated for increased precision. At each fold of the process, the error of the initial guess is halved, so that after folding through the procedure the algorithm produces a fold very close to the desired exact value. 

Given a fraction $0<s/t<1$ with odd denominator, the algorithm starts by folding an initial guess crease around $s/t$. Call this first crease's numerator $s_0=s$. For a previously formed crease $s_{n-1}$, the algorithm proceeds as follows:
\begin{enumerate}
\item If $s_{n-1}$ is odd, fold the right edge to the previous crease, forming a new crease with numerator $s_{n}=s_{n-1}+(t-s_{n-1})/2$.
\item If $s_{n-1}$ is even, fold the left edge to the previous crease, forming a new crease with numerator $s_{n}=s_{n-1}/2$.
\item Iterate the previous steps until for some $n$, $s_n=s_0$, or else repeat again the procedure if further precision is desired.
\end{enumerate}
That the procedure produces cycles can be seen as follows. Since $t$ is odd, either $s_{n-1}$ will be even or else $t-s_{n-1}$ will be even, and so at each step, a positive integer $s_n$ will be produced in the range $0<s_n<t$. After a maximum of $t$ steps, by the pigeonhole principle the $s_n$ will repeat a previous value, hence producing a cycle. 

We can analyze the algorithm and its cycles by rewriting the first two steps of the algorithm in the backwards direction as follows
\begin{enumerate}
\item $s_{n-1} = 2 s_n$, or
\item $s_{n-1} = 2 s_n - t$
\end{enumerate}
The backwards procedure is equivalent to multiplying $s_n$ by 2 mod $t$, that is, the procedure is an origami realization of cyclotomic cosets of 2 mod $t$. 

In the example of $\frac{1}{5}$ in Figure~\ref{figure2}, all possible creases are formed, that is, all numerators between 1 and $t$ are visited in the cycle. This isn't always the case, so that sometimes much shorter folding sequences can be found. For example, the cycle for $\frac{1}{7}$ is $\frac{1}{7}$, $\frac{4}{7}$, $\frac{2}{7}$, $\frac{1}{7}$, $\ldots$. Notice that in this case, the missing fractions form another independent cycle $\frac{3}{7}$, $\frac{5}{7}$ $\frac{6}{7}$, $\frac{3}{7}$,$\ldots$. 

The number of cycles that exist for a given denominator is given by the OEIS sequence A006694~\cite{oeisA006694}. The larger the number of cycles that exist for a given denominator, the smaller the cycle lengths on average, and hence the smaller a folding sequence. In fact, we will see below that for the fraction 4703/33215 only 36 folds are needed to form the wanted cycle, 4 orders of magnitude smaller than what one would naively expect, given the denominator.

As a means of describing the folding procedure, we will use R or L to denote folding the right or left edges of the paper to the previously folded crease, respectively. In the example of Figure~\ref{figure2} for finding 1/5, the sequence of folds from the initial guess is RRLL. Note that the sequence could have been stopped after just two steps (RR) if one is willing to rotate the strip of paper $180^{\circ}$ to find a new 1/5 approximant.

We seek to find a sequence of folds using Fujimoto's algorithm which will approximate $\pi$ using rational approximants to $\pi$. Given a choice of rational approximant, as long as enough folds or repetitions of the procedure are made, then the algorithm will coverge exponentially to the approximant such that the error in folding is smaller than the error coming from the rational approximation itself. The best rational approximants of $\pi$ are those coming from its continued fraction expansion, and we mostly consider these except for two cases below. 

It seems reasonable to only find $\pi$'s fractional part $0.14159\ldots$, since otherwise a much larger strip of paper is necessary, with the first 3 units of paper going unused. The first continued fraction convergent for $\pi$ is 22/7, the fractional part of which is given by 1/7, which has the simple sequence
\begin{equation}
\mathrm{RLL}
\end{equation}
This sequence can be repeated multiple times to decrease the error of the distance between the initial guess and the approximant. The approximant is accurate to 2 digits, $0.142\ldots$, and repeating this sequence twice will reduce any initial error by a factor of 64, which is likely sufficient precision.  

The next two smallest sequences which improve upon this are not produced by the continued fraction convergents, which, though they give more accurate results, have sequences which are longer and are presented afterwards. The next two sequences are for $18/127 = 0.1417\ldots$ with the 7 step sequence
\begin{equation}
\mathrm{LRLLRLL}
\end{equation}
and for $93/657 = 0.14155\ldots$ with the 18 step sequence
\begin{equation}
\mathrm{RRLLRRRRLLLLRLLRLL}
\end{equation}
Given the error of the initial guess fold, the first sequence of 7 steps will have decreased that error by a factor of $2^{-7}\approx 7.8\times10^{-3}$, which may not be sufficient for achieving the 3 digits of accuracy of the approximant. The sequence could be repeated a second time, for a total of 14 steps, so that the initial error is decreased by a factor of $2^{-14}\approx 6.1\times10^{-5}$ to ensure the accuracy of all digits. For the sequence with 18 steps, the initial error will be decreased by a factor of $2^{-18}\approx 3.8\times10^{-6}$, which will ensure that all 4 accurate digits of the approximant are achieved.

The next best sequence uses the continued fraction convergent 355/113. The fractional part is given by $16/113=0.1415929\ldots$ and from an initial guess we have the following sequence
\begin{equation}
\mathrm{LLLLRRRRLRRLRR}
\end{equation}
After this sequence of 14 folds, the last crease will be on the right as an approximant to $1-16/113$, so that the paper could be rotated $180^{\circ}$ to find an approximant for 16/113. Whatever the initial error of the guess, that error will have decreased by $2^{-14}\approx 6.1\times10^{-5}$, which may be sufficient for achieving 4 or 5 decimals of accuracy from the possible 6 decimals of the approximant, in which case this sequence of 14 steps is superior to the 18 step sequence above. If rotating the paper is undesirable, however, one could repeat the sequence in reverse, for a total of 28 steps, to arrive back on the left for a much more precise fold near the starting fold,
\begin{equation}
\mathrm{LLLLRRRRLRRLRR~~RRRRLLLLRLLRLL}
\end{equation}
Whatever the initial error of the guess, that error will have now decreased by $2^{-28}\approx 3.7\times10^{-9}$, so if the full accuracy of 6 digits of the convergent is desired, the sequence of length 28 can be used. 

Note that the full sequence of 28 folds contains a subsequence of 6 folds in succession all from the same side. In fact, four of the six happen immediately after the first half-sequence of 14 folds that end with a crease around 0.14159 from the right edge. The next four folds will subdivide this small length near the right edge of the paper by a factor of 16. Therefore, the paper used for this sequence should have a length greater than 11.3~cm, ideally at least double, so that none of the folds are within less than 1~mm from the edge of the paper. Standard Letter or A4 paper is almost triple this size along the long side. 

The next best folding sequence is for the convergent 104348/33215, accurate to 9 digits, which requires a sequence of length 36 to fold 4703/33215,
\begin{equation}
\mathrm{RRRRRLLRLLLRLRLRLRRLRRRRRRLLLLRLLRLL}
\end{equation}
In this sequence of folds, there is also a subsequence of 6 successive halvings, requiring a similarly sized length of paper compared to the previous sequence. 

Beyond these convergents, the lengths of the sequences of folds grow very quickly and there are greater subsequences of successive halvings, requiring ever larger lengths of paper to fold. 

In practicality, the precision of the half-sequence of 14 folds for the convergent 16/113 may be at the limit of what can be measured, since even with a 10 meter strip of paper, only 5 digits of precision can be measured if measuring to the nearest milimeter.

\section{Folding trigonometric approximants to \texorpdfstring{$\pi$}{Pi}}
Fujimoto noted that his algorithm also works for finding rational approximations to angles~\cite{fujimoto1983,fujimoto198X}. Given the right angle of a corner of a piece of paper, we would like to find a rational approximation to the angle $\theta\approx72.34321^{\circ}$ such that $\tan\theta=\pi$. Since $\theta/90^{\circ}=0.803813\ldots$, we can use the angular version of the algorithm to find rational approximants of this angle. 

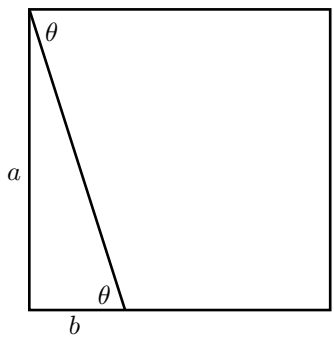
\begin{figure}[htpb]
\begin{center}
\begin{tikzpicture}[scale=1]
\draw[line width = 1pt] (0,0) rectangle (4,4);
\draw[line width = 1pt] (0,4) -- (1.2732,0);
\node at (.3,3.7) {$\theta$};
\node at (1.,0.2) {$\theta$};
\node at (-.2,1.8) {$a$};
\node at (0.6,-0.2) {$b$};
\end{tikzpicture}
\end{center}
\caption{Rational angle approximations can be found for the angle $\theta\approx 72.34321^{\circ}$ in the top left corner, for which $a/b=\pi$. \label{figure3}}
\end{figure}

We will assume the configuration shown in Figure~\ref{figure3}, and we will use the notation where L means the left edge of the paper is folded to the previous crease (such that the crease goes through the top left corner), while R means the top edge is folded to the previous crease (the crease going through the top left corner). Unfortunately, much higher accuracy is needed in the rational approximation to the angle in order to achieve a subsequent desired accuracy in the tangent of the angle.

The following sequence of 8 steps will find the rational angle approximation to $41/51=0.80392\ldots$
\begin{equation}
\mathrm{LRLLRRLL}
\end{equation}
which subsequently yield the approximation $3.1434\ldots$, accurate to 2 digits. The sequence of 8 steps will reduce any initial angle error by a factor of $2^{-8}\approx 3.9\times 10^{-3}$, which is only sufficient for achieving one digit of accuracy in the tangent of the angle. To achieve the full accuracy of the approximation for $\pi$, the sequence can be repeated, so that any angle error is reduced by a factor of $2^{-16}\approx 3.2\times 10^{-5}$.

The next best sequence is given by the rational angle approximation to $10034/12483 =0.80381318\ldots$, with sequence
\begin{equation}
\mathrm{RLRLLRRRLLLRLLRRLL}
\end{equation}
The rational angle approximation will yield the tangent approximation $3.141587\ldots$ with 4 digits of accuracy.

The next best sequence is for the rational angle approximation to $237047/294903=0.80381345\ldots$. It has 30 steps, but unfortunately it begins by bisecting the initial angle 9 times in a row, which is beyond the limits of practicality due to the difficulty of folding small angles. Laying aside this sequence, the next best sequence, which is a slightly worse approximation, is a sequence of 36 steps for the rational angle approximation to $80096/99645=0.8038135$, given by
\begin{equation}
\mathrm{RRRRRLRLLRLRLRRLLLRLRLLRRRLLLRLLRRLL}
\end{equation}
This gives the tangent approximation $3.1415937\ldots$ with 5 digits of accuracy, and the length of the sequence is sufficient to achieve this accuracy regardless of the initial error. Though easier than the rejected sequence, it still features 5 repeat bisections, which is difficult to fold in practice. 

Longer sequences for more accurate tangent approximations can be sought, but the length grows quickly and they may feature too many sequential bisections, making them impracticable to fold.

\section{Folding algebraic approximants to \texorpdfstring{$\pi$}{Pi}}
The final way we consider to fold an approximant of $\pi$ is to seek a polynomial equation of degree 2 or 3 with integer coefficients which has a real root close to $\pi$. The desired real root of the polynomial can then be found via the origami implementation of Lill's method, as detailed in~\cite{alperin2009}.

In practice, quadratics and cubics with ``small" integer coefficients are preferrable, e.g. coefficients with absolute value less than 64 or 128, since rational fractions of the side of the square of paper with denominator $2^n$ must be found, and divisions of the side of a square beyond 128 equal parts is not very convenient in practice. With these constraints, an exhaustive search for quadratics and cubics with integer coefficients $<100$ was performed. The best approximants come from cubics rather than quadratics, which is reasonable to expect given the larger number of cubic polynomials. Of those with coefficients $<64$, the best is
\beq
9x^3-9x^2-44x-52=0
\eeq
with its real root having a difference of $1.58\ldots\times 10^{-7}$ from $\pi$. 

Though the construction of the corresponding path of Lill's method can fit within a square with at most a division into 64 equal parts, carrying out the origami implementation of the method requires extra auxiliary crease lines which would fall outside of the bounds of the paper. The best cubic approximation whose origami construction falls within the square is
\beq
4x^3-22x^2+29x+2=0
\eeq
This cubic has a real root a distance of $3.42\ldots\times10^{-7}$ away from $\pi$, although it has another root $2.42\ldots$ which is close enough to the desired root that it could be mistakenly found rather than the $\pi$ approximant.

The best cubic approximant whose origami construction fits inside a square with at most 64 divisions and which cannot be mistaken for another nearby root is the cubic
\beq
8x^3-17x^2-16x-30=0
\eeq
which only has one real root, a distance of $1.27\ldots\times10^{-5}$ from $\pi$. 

Finally, the best cubic approximation of $\pi$ was found which can fit within a Letter or A4 sheet of paper such that after folding the approximant the root can be measured directly in cm is the cubic given by  
\beq
x^3-x^2-8x+4=0
\eeq
which has the real root $3.1413\ldots$, accurate to 3 digits.

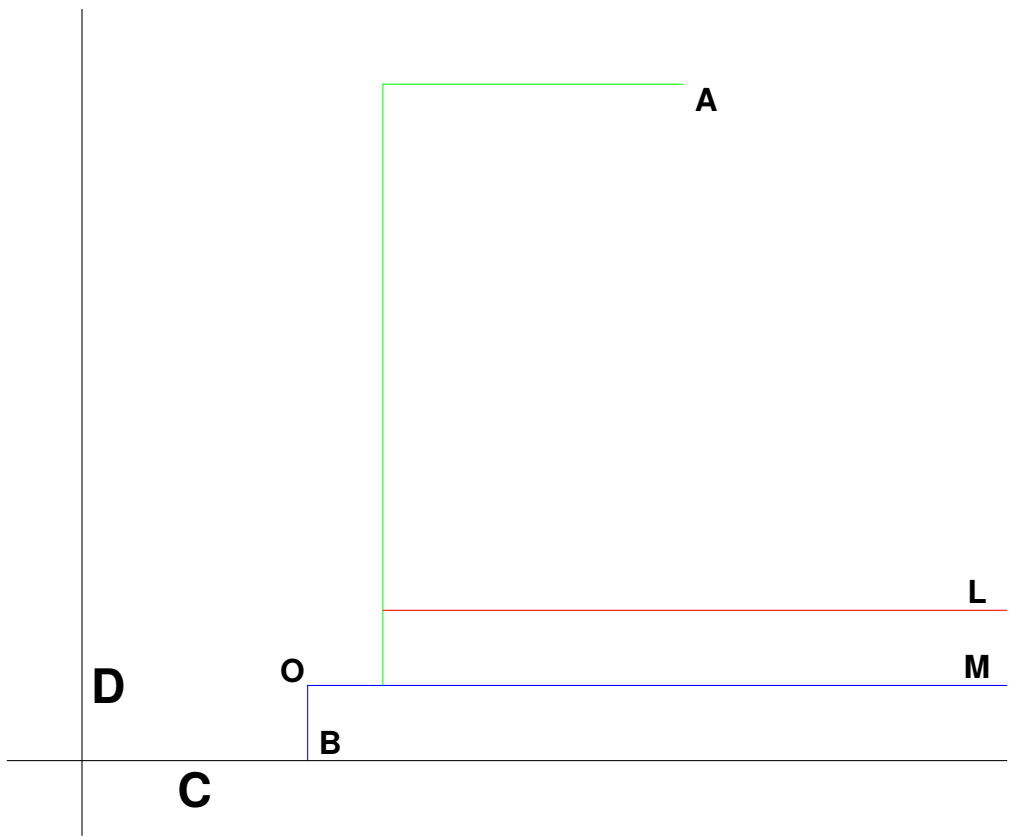
\begin{figure}[htpb]
\begin{center}
\begin{tikzpicture}[scale=0.4]{x=1cm,y=1cm}
\begin{scope}[shift={(4.5,-8.5)}]
\draw[blue,line width=0.4pt] (0,0) -- (0,1) -- (9.3,1);
\draw[green,line width=0.4pt] (1,1) -- (1,9) -- (5,9);
\draw[red,line width=0.4pt] (1,2) -- (9.3,2);
\draw[black,line width=0.4pt] (-4.,0) -- (9.3,0);
\draw[black,line width=0.4pt] (-3,-1) -- (-3,10);
\node at (-0.25,1.25) {\fontfamily{qhv}\selectfont\textbf{\large{O}}};
\node at (5.3,8.9) {\fontfamily{qhv}\selectfont\textbf{\large{a}}};
\node at (0.35,0.35) {\fontfamily{qhv}\selectfont\textbf{\large{b}}};
\node at (8.9,2.35) {\fontfamily{qhv}\selectfont\textbf{\large{L}}};
\node at (8.9,1.35) {\fontfamily{qhv}\selectfont\textbf{\large{M}}};
\node at (-1.5,-0.55) {\scalebox{1.5}{\fontfamily{qhv}\selectfont\textbf{\large{C}}}};
\node at (-2.55,1.) {\scalebox{1.5}{\fontfamily{qhv}\selectfont\textbf{\large{D}}}};
\end{scope}
\end{tikzpicture}
\end{center}
\caption{Lill's method for finding the real root $3.1413\ldots$ of the cubic $x^3-x^2-8x+4=0$. A single fold is made by aligning point a to line D while at the same time aligning point b to line L, making a crease that crosses line M. The measurement in cm from O to M will give the $\pi$ approximant when the length from b to O is 1~cm. \label{fig:foldingactivity}}
\end{figure}
The origami folding procedure is given in Figure~\ref{fig:foldingactivity}, where O is the origin. The point b is at (0,-1), the line segment L is at y=1, the point a is at (4,8), and the line D is at x=-4. A single fold is made by aligning point a to line D while at the same time aligning point b to line L, making a crease that crosses line M. A printing such that the length from b to O is 1~cm will result in a crease crossing line M such that the length from O to the crease intersection is $3.1413\ldots$. Note that in practice care must be taken to ensure that a printer will print at 100\% scaling.

\section{Folding \texorpdfstring{$\pi$}{Pi} exactly using curved creases}
\label{sec:piexact}
\begin{figure}[htpb]
\begin{center}
\includegraphics[height=2.in]{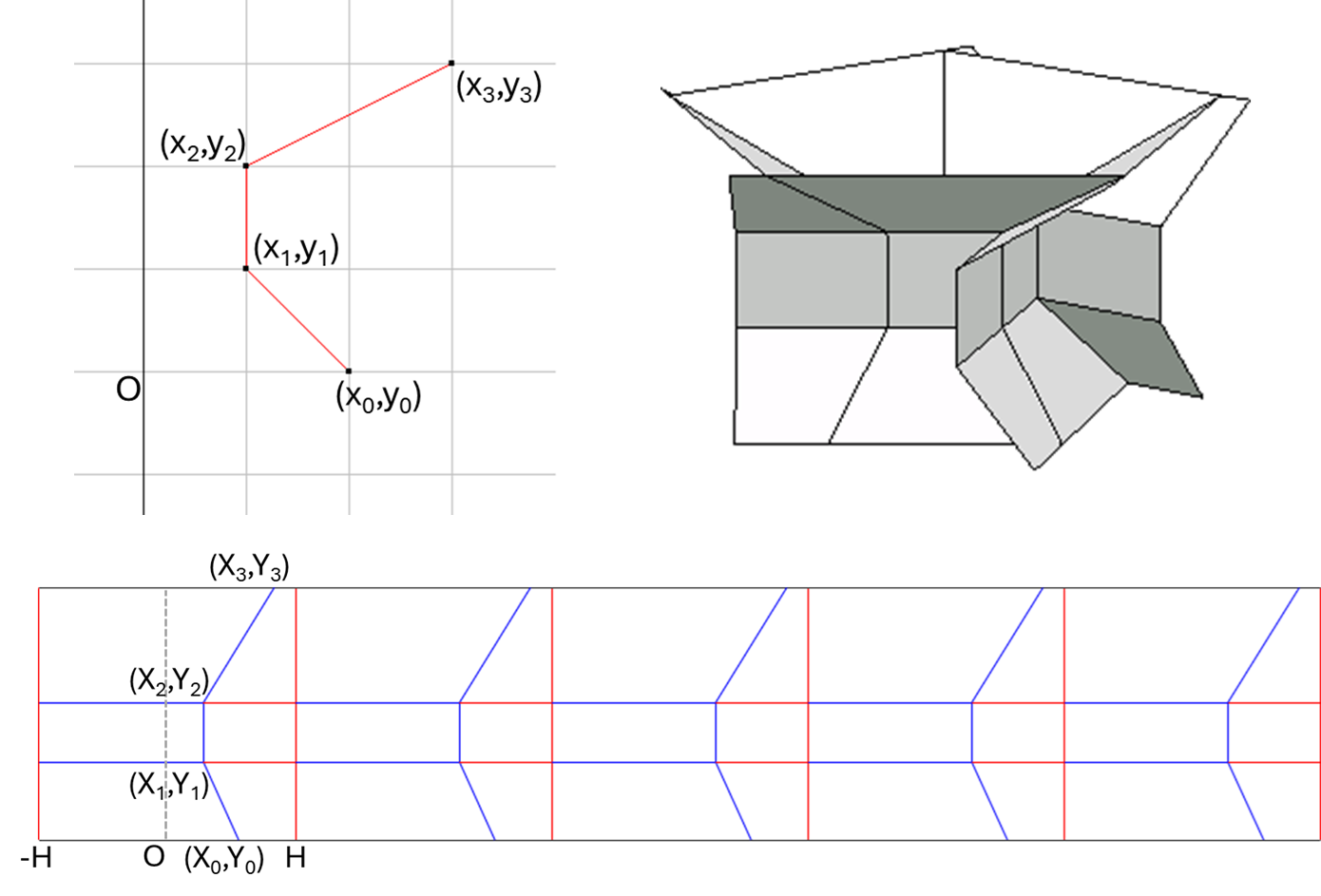}
\end{center}
\caption{Annotated cropped windows of Jun Mitani's Ori-Revo software showing the silhouette (top left), folded model (top right), and crease pattern (bottom).\label{fig:orirevo}}
\end{figure}

Jun Mitani has popularized rotationally symmetric origami models through his software Ori-Revo~\cite{mitani2024}, a sample of which is shown in Figure~\ref{fig:orirevo}. In the software, the user can select points, connected by lines, which represent the silhouette of the shape after folding along a crease pattern that is produced. The number of segments and other parameters can be chosen by the user. The mathematics involved in producing the crease pattern is described in detail in Robert Lang's book Twists, Tilings, and Tessellations~\cite{lang2018}, of which we focus here on the cylindrical thin-flange case. 

A limitation of the software and the mathematical description in Twists, Tilings, and Tessellations, however, is that only line segments are allowed or presented. We can easily adapt the equations to the continuous case, however. Given a silhouette curve $x(y)$, the crease pattern curve can be calculated as:
\beqr
X &=& x(y)\tan\left(\frac{\pi}{n}\right)\\
Y &=& \int_{0}^{y}\sqrt{1+\left(\frac{dx}{dy} \right)^2}dy \\
H &\geq& \max(X)
\eeqr
where $X$ and $Y$ are crease pattern coordinates given parametrically, and where $H$ is the half-width of each panel of the crease pattern segments. The greater the quantity $H$, the larger the flange that sticks out of the folded shape, with the equality meaning that the flange will disappear for at least one point along the silhouette curve. We don't need to restrict our original curve to be a function, however. We can also calculate a crease pattern for a parametrically defined curve given by $(x(t),y(t))$ as
\beqr
X &=& x(t)\tan\left(\frac{\pi}{n}\right)\\
Y &=& \int_{0}^{t}\sqrt{\left(\frac{dx}{dt} \right)^2+\left(\frac{dy}{dt} \right)^2}dt \\
H &\geq& \max\left(X \right)
\eeqr
The author has used these equations to produce curved crease patterns in Mathematica~\cite{mathematica2024}.  

The formula for $Y$ in both cases is in fact just the arclength of the curve, and in certain cases the formula can produce closed-form solutions. These solutions are in general transcendental, so that constructing the crease pattern curve to fold is itself non-trivial. 

However, we can work the other way around, asking what silhouette is produced by folding along a given chosen curve. Given a crease pattern function $X(Y)$, with a chosen $H \geq \max(X)$, the silhouette equations become,
\beqr
x &=& X(Y)\cot\left(\frac{\pi}{n}\right)\\
y &=& \int_{0}^{Y}\sqrt{1-\left(\frac{dX}{dY} \right)^2}dY 
\eeqr
If the crease pattern curve is given parametrically as $(X(t),Y(t))$, the silhouette equations become
\beqr
x &=& X(t)\cot\left(\frac{\pi}{n}\right)\\
y &=& \int_{0}^{t}\sqrt{\left(\frac{dY}{dt} \right)^2 - \left(\frac{dX}{dt} \right)^2}dt 
\eeqr

We can now construct the curves we wish and ask what final folded shapes result. One of the simplest curves other than line segments is a parabola. For the parabola $X=Y^2/2$ between $0\leq Y\leq1$, the resulting folded silhouette is given by the curve
\beqr
x &=& \frac{Y^2}{2}\cot\left(\frac{\pi}{n}\right)\\
y &=& \frac{1}{2}\left[Y\sqrt{1-Y^2} + \sin^{-1}(Y) \right]
\eeqr
For convenience, we can choose $n=4$, so that 
\beqr
x &=& \frac{Y^2}{2}\\
y &=& \frac{1}{2}\left[Y\sqrt{1-Y^2} + \sin^{-1}(Y) \right]
\eeqr
or equivalently,
\beq
y = \frac{1}{2}\left[\sqrt{2x}\sqrt{1-2x} + \sin^{-1}(\sqrt{2x}) \right]
\eeq

\begin{figure}[htpb]
\begin{center}
\subfloat[]{
	\centering
	\begin{tikzpicture}[scale=0.7]
	\begin{axis}[
		axis equal image,
	    xmin=0,xmax=1,axis x line=middle,
	    ymin=0,ymax=1.,axis y line=center,
	     ] 
	     \addplot [domain=0:1,samples=100,color=red]({1/2*x^2},{x});
	 \end{axis}
	\end{tikzpicture}
}
\hspace{0.04\textwidth}
\subfloat[]{
	\centering
	\begin{tikzpicture}[scale=0.7]
	\begin{axis}[
		axis equal image,
	    xmin=0,xmax=1,axis x line=middle,
	    ymin=0,ymax=1.,axis y line=center,
	     ] 
	     \addplot [domain=0:1,samples=100,color=red]({x^2/2},{1/2*(x*sqrt(1-x^2)+asin(x)/180*pi)});
	     \draw[color=blue,dashed] (0,pi/4) -- (0.5,pi/4);
		 \draw[color=blue,dashed] (0.5,0) -- (0.5,pi/4);
		\node at (0.6,pi/4+0.05) {$(0.5,\pi/4)$};
		\node[circle,fill,inner sep=2pt] at (0.5,pi/4) {};
	 \end{axis}
	\end{tikzpicture}
}
\end{center}
\caption{On the left the plot of the crease pattern function $x=y^2/2$ between $0\leq x \leq 1$. On the right the plot of the silhouette function $y=1/2\left[\sqrt{2x}\sqrt{1-2x} + \sin^{-1}(\sqrt{2x}) \right]$ between $0\leq x \leq 0.5$.\label{fig:picurves}}
\end{figure}
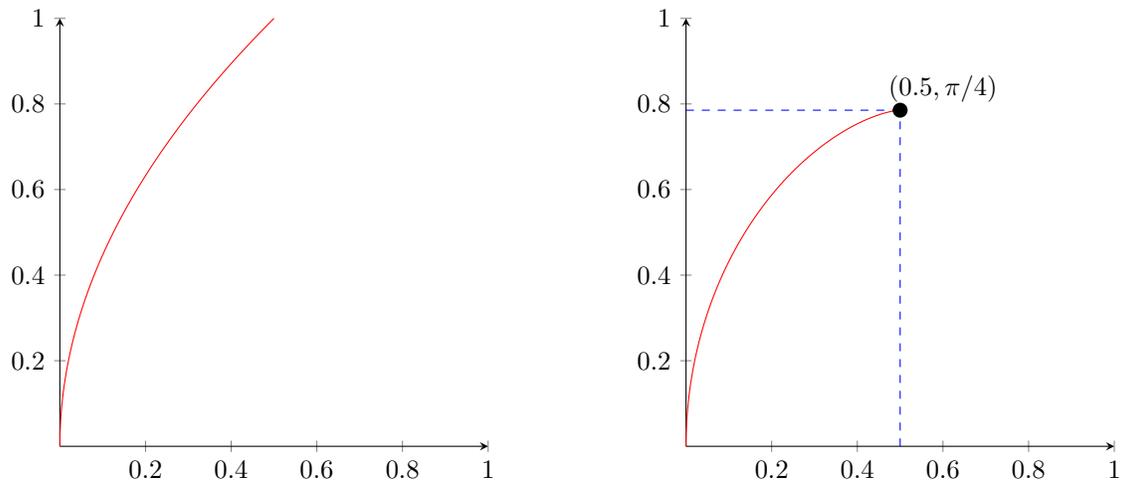

\begin{figure}[htpb]
\begin{center}
\includegraphics[height=2.25in]{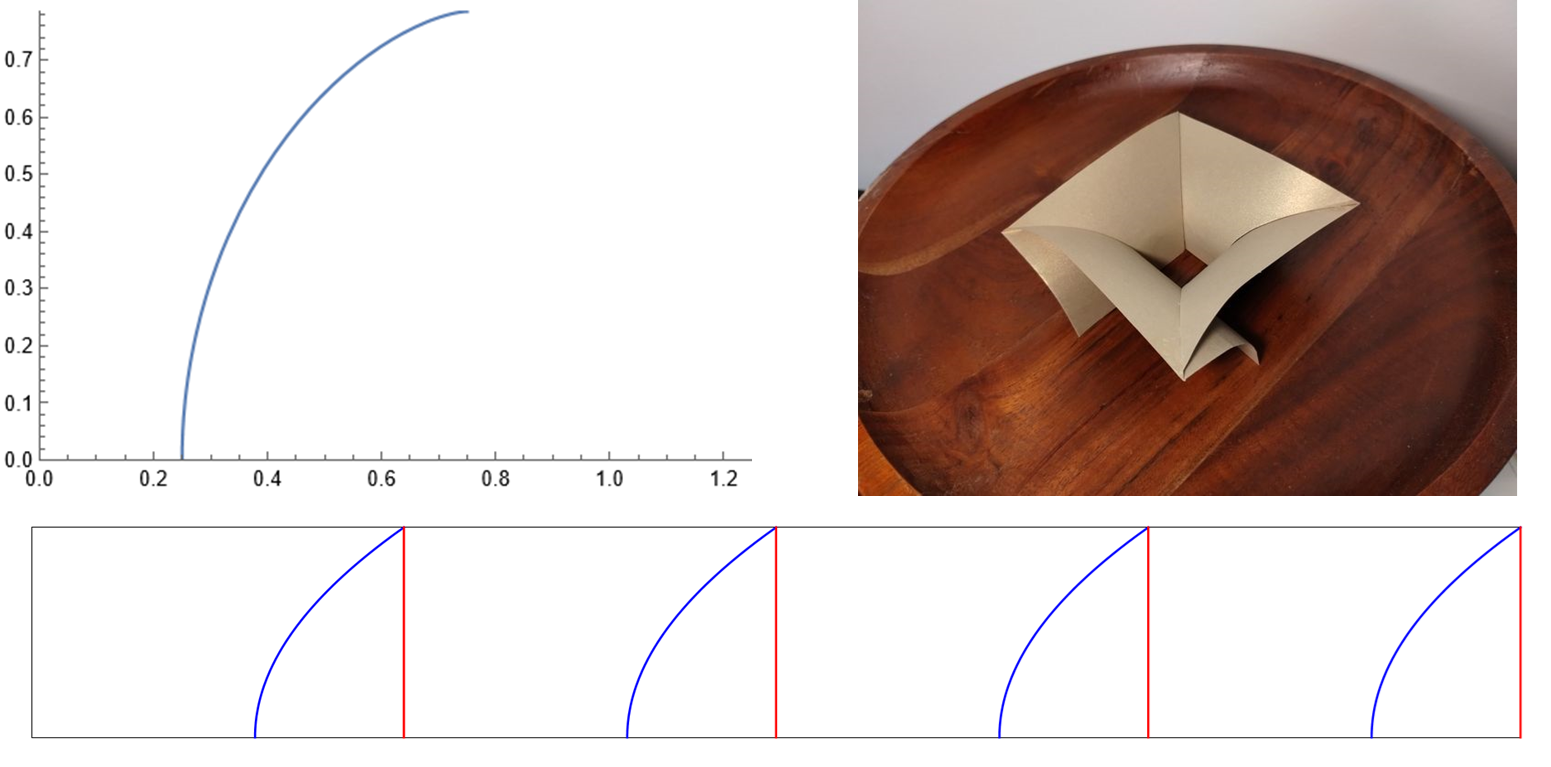}
\end{center}
\caption{Silhouette design of parabola $x=\frac{y^2}{2}$ (top left) with 4 segments producing $\pi$ after folding, crease pattern (bottom), folded model by the author (top right).\label{fig:piquad}}
\end{figure}

We can now consider the endpoint coordinate $Y=1$, at $x=1/2$, and we find that 
\beq
y\left(\frac{1}{2}\right) = \frac{\pi}{4}
\eeq
The crease pattern curve and the folded silhouette curve are shown in Figure~\ref{fig:picurves}. The full crease pattern is shown in Figure~\ref{fig:piquad} with the author's folding of the crease pattern. 

Since we know the crease pattern being folded must have a height of 1 before folding in the crease pattern coordinates, we can simply find the ratio of the height after folding to the height before folding, with the result being $\pi/4$. And since we will have folded four segments to achieve the result, the construction given by the crease pattern shown in Figure~\ref{fig:piquad} in fact folds $\pi$! And because the surface rulings are parallel to the ground and the cross-sections are therefore squares, one might say, tongue-in-cheek, that ``origami can square the circle". 

\section{Other transcendental numbers}
The number $\pi$ is not the only transcendental number we can consider using the methods above. We can seek to approximate any other number using continued fraction convergents, or trigonometric or algebraic approximants. More interesting is to consider other convergent sequences, but the author is unaware of other constructions which could be adapted to origami that are similar to the one above for $\pi$. 

The curved crease method to produce $\pi$ can also produce other transcendental numbers. For example, folding along the cubic $X=Y^3/3$ between $0\leq Y \leq 1$ for 4 segments produces the transcendental number
\beq
y(Y=1) = \frac{\Gamma\left(\frac{1}{4} \right)^2}{6\sqrt{2\pi}} \label{cubicgamma}
\eeq
and we note that the expression (\ref{cubicgamma}) with a different rational constant, 1/8 instead of 1/6, was previously noted in connection with an infinite product associated with the regular paperfolding sequence, also known as the dragon curve sequence~\cite{allouche2015}. 

In general we have the following silhouette equations for folding with 4 segments,
\beqr
X &=& \frac{Y^n}{n}\\
y &=& \int_{0}^{Y}\sqrt{1+Y^{2n-2}}dY \\
&=& {}_2F_1\left(-\frac{1}{2},\frac{1}{2n-2};1+\frac{1}{2n-2};Y^{2n-2} \right) \\
x &=& \frac{Y^n}{n}
\eeqr
where ${}_2F_1(a,b;c;x)$ is the Gauss hypergeometric function. 

The endpoint at $Y=1$ is given by
\beq
y(Y=1) = \frac{\sqrt{\pi}\Gamma\left(\frac{1}{2n-2}\right)}{2n\Gamma\left(\frac{1}{2}+\frac{1}{2n-2}\right)}= \frac{2^{\frac{1}{n-1}}\Gamma\left(\frac{1}{2n-2}\right)^2}{4n\Gamma\left(\frac{1}{n-1}\right)}
\eeq
where $\Gamma(x)$ is the Gamma function. 

This procedure has the potential to construct other transcendental constants, although it requires first constructing polynomials of high degree. Not much is known about the algebraic independence of the numbers $\Gamma\left(1/k\right)$ except that the numbers for $k=2,3,4$ are known to be algebraically independent~\cite{waldschmidt1977}. This means that the cubic procedure which produces a ratio involving $\Gamma(1/4)$ and $\Gamma(1/2)=\sqrt{\pi}$, actually constructs a new transcendental constant $\Gamma(1/4)$ distinct from $\pi$. 

\section{Conclusions}
We have presented several methods for constructing $\pi$ either exactly, convergently, or approximately. Of course, all of these methods have the inherent limitation of precision of both folding and measuring, let alone material properties. In addition, constructing a continuous parabola to fold along in order to construct $\pi$ exactly itself falls outside of the realm of origami. 

This isn't the first time that $\pi$ has been considered in origami. In~\cite{hull2007}, a semi-circle was folded in a strip of paper and then the edge of the paper was wrapped around the semi-circle in order to measure its perimeter. This procedure requires very careful alignment of the edge with the curved crease while the entire paper is folded in 3D, and it would be equivalent to simply rolling paper around a cylinder to measure its perimeter, although everything occurs within the one sheet of paper itself. By way of contrast, our procedure does not involve any alignments, although a cylindrical sheet of paper is required. The simple act of folding deforms the paper so that its resultant height is $\pi/4$. Therefore, we consider our procedure a new non-trivial folding of $\pi$, due to the the paper exercising its intrinsic geometry under the curved crease of a parabola to form $\pi$. 

Of course, because the Huzita-Justin straight crease axioms can construct solutions to quadratics, they can be used in principle to construct a dense set of origami-constructible points along the parabola. Therefore, constraining ourselves to only points constructible using standard origami procedures, we can define a new convergent sequence which approximates $\pi$, where more and more points along the parabola $x=y^2/2$ are constructed, and the height measured after folding along all of them will approach $\pi/4$. Further, since standard stright-line origami operations can also construct a dense set of points along cubics, the curved crease procedure also produces a convergent sequence for constructing the transcendental number $\Gamma\left(\frac{1}{4} \right)$. In fact, for any given curve having a dense number of origami-constructible points, an origami sequence can be defined which may converge to new transcendental numbers.

It remains, therefore, that from the perspective of straight-line crease origami, $\pi$ and other transcendental numbers retain their enchanting mystery as numbers that can only be produced at best convergently in a finite number of steps.

\section{Acknowledgements}
Most of the work in this paper was produced while the author was in the School of Information and Physical Sciences at the University of Newcastle, NSW, Australia in 2018 and the author wishes to thank the generous hospitality of Judy-anne Osborn and Richard Brent during this time.  


\end{document}